\documentclass[11pt,a4paper]{article}
\usepackage[latin1]{inputenc}
\usepackage{amsmath}
\usepackage{amsfonts}
\usepackage{amssymb}
\usepackage{graphicx}
\usepackage{graphics}
\usepackage[top=2cm, bottom=2cm, left=2cm, right=2cm]{geometry}
\usepackage{enumerate}
\begin{document}
\thispagestyle{empty}

\textbf{{\Large Spectral Methods And Primes Counting Problems}}\\
N. A. Carella\\
\vskip .5 in
\textbf{\textit{Abstract}:} A recent heuristic argument based on basic concepts in spectral analysis showed that the twin prime conjecture and a few other related primes counting problems are valid. A rigorous version of the spectral method heuristic, and a proof of the more general dePolignac conjecture on the existence of infinitely many primes pairs $p$ and $p + 2k$, $k \geq 1$ are considered within. Further, a new application of the spectral method to the existence of infinitely many quadratic twin primes $n^{2}+1$ and $n^{2}+3$, $n \geq 1$, is also proposed in this note.\\

\vskip .5 in
\tableofcontents

\vskip .5 in

\\

\pagenumbering{Page gobble}
\pagenumbering{arabic}

\newpage
\section{Introduction}
A heuristic argument in favor of the existence of infinitely many twin primes numbers $p$ and $q = p + 2$ was devised in \cite{GP99}, and generalized to some related linear prime pairs $p$ and $ap + bq = c$ counting problems, where $a, b, c \in \mathbb{Z}$ are fixed integers in \cite{GP06}. The heuristic argument is based on basic spectral techniques, as Ramanujan series, and the Wiener-Khintchine formula, applied to the correlation function $\sum_{n \leq x} \Lambda(n)\Lambda(n+2k)$. Extensive information on the twin primes problem appears in \cite{RP96}, \cite{NW00}, \cite{GG05}, \cite{GD09}, et cetera. 
\\

This work proposes rigorous versions of the formal proofs of several primes pairs counting problems, and related prime Diophantine equations such as $ax + by + c$ and $ax^2+bx+c$ over the set of integers $\mathbb{Z}$. The basic idea combines elements of spectral analysis and complex functions theory to derive a rigorous proof. The result for linear primes pairs is as follows.
\\

\textbf{Theorem 1.1.}  \textit{Let $k\geq 1$ be a fixed integer, and let $x\geq 1$ be a large real number. Then}
\\
\begin{equation}
  \sum_{n\leq x}\Lambda(n)\Lambda(n+2k) \geq  \alpha_{2} x+O\left (x \frac{(\log \log x)^2}{\log x} \right ),
\end{equation}
\\
\textit{for some $\alpha_{k} > 1/2$ constant.}                                   
\\

This implies that there is an effective asymptotic formula for the primes pairs counting function
\\
\begin{equation}
\pi_{2k}(x) =\#\{p,p+2k\leq x \}\geq\frac{1}{2} \frac{x}{\log^2 x}+O\left (x \frac{(\log \log x)^2}{\log^3 x} \right ).
\end{equation} \\

The subsequent Sections collect the necessary elementary results in arithmetic functions, spectral analysis, complex function theory, and some ideas on convergent series. Section 8 assembles a short proof of Theorem 1.1.
\\

The quadratic twin primes $n^2+1, n^2+3, n\geq 1,$ counting problem claims the existence of infinitely many such prime pairs, the heuristic appears \cite{HL23}. As in the previous case for linear twin prime pairs, the basic idea of the proof combines elements of spectral analysis and complex functions theory to derive a rigorous proof applied to the correlation function $\sum_{n^{2}+1\leq x}\Lambda(n^{2}+1)\Lambda(n^{2}+3)$, see \cite[p. 62]{HL23}, \cite[p. 343]{NW00}, \cite[p. 406]{RP96}, et alii. The result is as follows.\\

\textbf{Theorem 1.2.}  \textit{Let $x\geq 1$ be a large real number. Then}\\
\begin{equation}
\sum_{n^{2}+1\leq x}\Lambda(n^{2}+1)\Lambda(n^{2}+3) \geq  \alpha_{2} x^{1/2}+O\left (x^{1/2} \frac{(\log \log x)^2}{\log x} \right ),
\end{equation}\\
\textit{for some $\alpha_{2} > 1/2$ constant.}      \\

This implies that there is an effective asymptotic formula for the quadratic twin primes counting function\\
\begin{equation}
\pi_{f_2}(x) =\#\{n^{2}+1,n^{2}+3\leq x \}\geq\frac{1}{2} \frac{x^{1/2}}{\log^2 x}+O\left (x^{1/2} \frac{(\log \log x)^2}{\log^3 x} \right ).
\end{equation}\\
Section 9 assembles a short proof of Theorem 1.2.

\vskip 1 in
\section{Properties of the Ramanujan Sum and Identities}
The symbols $\mathbb{N} = \{ 0, 1, 2, 3, \ldots \}$ and $\mathbb{Z} = \{ \ldots, -3, -2, -1, 0,  1, 2, 3, \ldots \}$ denote the subsets of integers. For $n \geq 1$, the Mobius function is defined by
\\
\begin{equation}
\mu(n) =
\left\{
\begin{array}{ll}
(-1)^v     &n=p_{1}\cdot p_{2} \cdot \cdot \cdot p_{v},\\
0          &n\ne p_{1} \cdot p_{2}\cdot \cdot \cdot p_{v},\\
\end{array}
\right.
\end{equation}
\\  
where the $p_{i}$ are primes. The Euler totient function is defined by\\ 
\begin{equation}
\varphi(n)=\prod_{p|n}(1-1/p).
\end{equation} \\
And the Ramanujan sum is defined by\\ 
\begin{equation}
c_{q}(n)=\sum_{\gcd(k,q)=1}e^{i 2 \pi k n /q}.
\end{equation} \\
The basic properties of $c_{q}(n)$ are discussed in \cite[p. 160]{AP76}, \cite[p. 139]{HR59}, $\cite{HW08}$, et alii.

\subsection{Some Properties}
The arithmetic function $c_{q}(n)$ is  multiplicative but not completely multiplicative. Some basic properties of the function $c_{q}(n)$ are listed here.\\

\textbf{Lemma 2.1.}   (Identities)   \textit{ Let $n$ and $q$ be integers. Then }
\begin{enumerate}[(i)]
\item $c_{q}(n)=\mu(q/d) \varphi(q)/\varphi(q/d)$, \textit{where} $d = \gcd(n, q)$.
\item $c_{q}(n)=\sum_{d|\gcd(n,q)} \mu(q/d)d$.
\item $c_{q}(n)=\prod_{p^{\alpha}||n}c_{p^{\alpha}(n)}$,  \textit{where} $p^{\alpha}||n $ \textit{denotes the maximal prime power divisor of} $n$.
\item $c_{2m}(n)=-c_{m}(n)$ \textit{for any odd integers} $m, n \geq 1$.
\item $c_{2^{v}m}(n)=0$ \textit{for any odd integers} $m, n \geq 1$, \textit{and} $v \geq 2$.
\item \textit{For any prime power}
\begin{equation}
c_{p^v}(n)=\left \{\begin{array}{ll}
p^{v-1}(p-1) & \text{if } \gcd (n,p^v)=p^v, v\geq 1,\\
-p^{v-1} & \text{if } \gcd (n,p^v)=p^{v-1}, v\geq 1,\\
0 & \text{if } \gcd (n,p^v)=p^u \text{ with }0 \leq u\leq v-2. 
\end{array}\right.
\end{equation} 	
\end{enumerate}

\textbf{Proof:} (ii) implies (iv) and (v). As the integer $n$ is odd, the index of the finite sum $d | \gcd(2^{v}m, n)$ is odd. Thus, \\
\begin{equation}
 c_{2^{v}m}(n)=\sum_{d|\gcd(2^{v}m,n)} \mu(2^{v}m/d)d=\mu(2^{v})\sum_{d|\gcd(m,n)}\mu(m/d)d =\left\{
\begin{array}{ll}
-c_{m}(n)     &\text{if }v=1,\\
0          &\text{if }v\geq2,\\
\end{array}
\right.
\end{equation}
This completes the verifications of (iv) and (v).  $\blacksquare$\\

\textbf{Lemma 2.2.}   (i) \textit{For any fixed $n \geq 1$, the estimate $1 \leq | c_{q}(n)|\leq 2n \log \log n$  holds all $q \geq 1$.}
(ii) \textit{For any fixed $q \geq 1$, the estimate $1 \leq | c_{q}(n)|\leq 2q \log \log q$  holds all $n \geq 1$.}\\

\textbf{Proof:} Using the formulas $c_{q}(n)=\sum_{d|\gcd(n,q)} \mu(q/d)d$, and $\sigma(n)=\sum_{d|n}d \leq 2n \log \log n$ it is easy to confirm that \\
\begin{equation}
1 \leq | c_{q}(n)|= \left |\sum_{d|\gcd(n,q)} \mu(q/d)d \right  |\leq \sum_{d|n}d  \leq 2n \log \log n .
\end{equation} \\
This verifies the inequality (i).  $\blacksquare$\\

The detail on the upper bound $\sigma(n)=\sum_{d|n}d \leq 2n \log \log n$ appears in \cite[p. 116]{TG15} or similar references. A sharper bound  
$1 \leq | c_{q}(n)|\leq d(n)$ is claimed in \cite[p. 139]{HR59}. Here $d(n)=\sum_{d|n}1 \ll n^\varepsilon$, where $\varepsilon>0$ is an 
arbitrarily small number. \\

These and many other basic properties of $c_{q}(n)$ are discussed in \cite[p. 160]{AP76}, \cite[p. 237]{HW08}, etc. New applications and new proofs were developed in \cite{FK12} for various exponential finite sums. The followings are proved ibidem. \\

\textbf{Lemma 2.3.}  (\cite[Corollary 3.6]{FK12}).  \textit{If $n \in \mathbb{N}$ is an integer, then} \\
\begin{equation}
\sum_{d|n}c_d(n/d) =\left \{
\begin{array}{ll}
\sqrt{n}     &\text{if  }n=m^2, \\
0          &\text{if  }n\ne m^2. 
\end{array}
\right .
\end{equation} \\

\textbf{Lemma 2.4.}  (\cite[Corollary 3.12]{FK12}).  \textit{If $m \in \mathbb{N}$ is an integer, then} \\
\begin{equation}
\sum_{d|n}c_d(m) =\left \{
\begin{array}{ll}
n     &\text{if  }n|m, \\
0          &\text{if  }n\nmid m. 
\end{array}
\right .
\end{equation} \\

\vskip 1 in
\section{Mean Values of Functions}
A few definitions and results on the spectral analysis of arithmetic functions are discussed in this section. \\

Let $f : \mathbb{N} \longrightarrow \mathbb{C}$ be a complex valued function on the set of nonnegative integers. The mean value of a function is denoted by\\
\begin{equation}
M(f)=\lim_{x \rightarrow \infty} \frac{1}{x} \sum_{n \leq x}f(n),
\end{equation} \\
confer \cite[p. 46]{SS94}. \\

For $q \geq 0$, the \textit{q}th coefficient of the Ramanujan series $f(n)=\sum_{n \geq 1}a_q c_q(n)$ is defined by \\
\begin{equation}
a_q = \frac{1}{\varphi(q)} M(c_q(n) \overline{f}),
\end{equation} \\
where $\overline{f}$ is the complex conjugate of $f$, and $\varphi$ is the Euler totient function. \\

\textbf{Lemma 3.1.}  (Orthogonal relation).  \textit{Let $n, p$ and $ r \in \mathbb{N}$ be integers, and let $c_q(n)=\sum_{\gcd(k,q)=1}e^{i 2 \pi k n /q}$.} 
\textit{Then} \\
\begin{equation}
M(c_q(n) c_r(n))=\left \{
\begin{array}{ll}
\varphi(q)     &\text{if  }q=r, \\
0          &\text{if  }q\ne r. 
\end{array}
\right .
\end{equation} \\
\textbf{Proof:} Compute the mean value of the product $c_q(n)c_r(n)$:
\begin{eqnarray}
M(c_q(n)c_r(n))&=&\lim_{x \rightarrow \infty} \frac{1}{x} \sum_{n \leq x}c_q(n)c_r(n)\\
&=&\lim_{x \rightarrow \infty} \frac{1}{x} \sum_{n \leq x,}\sum_{\gcd(u,q)=1}e^{i 2 \pi  n \frac{u}{q}} \sum_{\gcd(v,r)=1}e^{i 2 \pi  n \frac{v}{r}}\\
&=&\lim_{x \rightarrow \infty} \frac{1}{x} \sum_{n \leq x,} \sum_{\gcd(u,q)=1,}\sum_{\gcd(v,r)=1}e^{i 2 \pi  n(\frac{u}{q}- \frac{v}{r})}.
\end{eqnarray} \\

If $u/q=v/r$, then the double inner finite sum collapses to $\varphi(q)$. Otherwise, it vanishes:\\
\begin{eqnarray*}
\lim_{x \rightarrow \infty} \frac{1}{x} \sum_{n \leq x,} \sum_{\gcd(u,q)=1,}\sum_{\gcd(v,r)=1}e^{i 2 \pi  n(\frac{u}{q}-\frac{v}{r})}
&=& \sum_{\gcd(u,q)=1,}\sum_{\gcd(v,r)=1}\lim_{x \rightarrow \infty}  \sum_{n \leq x}\frac{1}{x}e^{i 2 \pi  n(\frac{u}{q}-\frac{v}{r})}\\
&=& \sum_{\gcd(u,q)=1}\sum_{\gcd(v,q)=1}\lim_{x \rightarrow \infty}  \frac{1}{x} \left ( \frac{1-e^{i 2 \pi  (\frac{u}{q}-\frac{v}{r})(x+1)}}{1-e^{i 2 \pi  (\frac{u}{q}-\frac{v}{r})}} \right ) \\
&=&0.
\end{eqnarray*} 
\begin{equation}
\quad
\end{equation}\\
for $u/q \ne v/r$.   $\blacksquare$ \\

Similar discussions are given in \cite[p. 270]{SS94}, and the literature. Extensive details of the analysis of arithmetic functions are given 
in \cite{HR59}, \cite{SS94}, \cite{LL10}, \cite{RM13} and the literature.

\vskip 1 in
\section{Some Harmonic Series For Arithmetic Functions}
The vonMangoldt function is defined by\\
\begin{equation}
\Lambda(n) =
\left\{
\begin{array}{ll}
\log p&n=p^k,\\
0&n\ne p^k,k\ge1.\\
\end{array}
\right.
\end{equation} \\

There are various techniques for constructing analytic approximations of the vonMangoldt function. The Ramanujan series and the Taylor series considered here are equivalent approximations. \\

\subsection{Ramanujan Series}
The Ramanujan series $f(n,s)$ at $s=1$ coincides with the ratio $\frac{\varphi(n)}{n}\Lambda(n)=\Lambda(n)+O(\frac{\log n}{n})$. Ergo, this series realizes an effective analytic approximation of the vonMangoldt function $\Lambda(n)$ as a function of two variables $n\geq2$ and $s\approx 1$. \\

The convergence property of the kernel series  $f(n,1)$ was proved in an earlier work, it should be noted that it is not a trivial calculation. \\

\textbf{Lemma 4.1.}   (\cite[p. 25]{KL12})  \textit{The Ramanujan series }\\
\begin{equation}
\sum_{q\geq 1}\frac{\mu(q)}{ \varphi(q)}c_{q}(n)
\end{equation} \\
\textit{diverges for $n = 1$ and conditionally converges for $n \geq 2$.} \\

\textbf{Lemma 4.2.} \textit{Let $n \geq 2$ be a fixed integer. For any complex number $ s \in \mathbb{C}, \Re e(s)>1 $, the Ramanujan series}\\
\begin{equation}
f(n,s)= \sum_{q\geq 1}\frac{\mu(q)}{q^{s-1}\varphi(q)}c_{q}(n)
\end{equation} \\
\textit{is absolutely and uniformly convergent.}\\

\textbf{Proof:} Suppose \(n\geq 1\) is a fixed integer. The inequality in Lemma 2.2 gives\\ 
\begin{equation}
\left | \sum _{q\geq 1} \frac{\mu (q)}{q^{s-1}\varphi(q)}c_q(n) \right| \leq \sum _{q\geq 1} \frac{1}{q^{\sigma -1}\varphi (q)} 
\left| c_q(n) \right| \leq (2n \log \log n) \sum _{q\geq 1} \frac{1}{q^{\sigma -1}\varphi (q)}. \end{equation}\\
Now, using \(q/ \log q \leq \varphi(n)\leq q \), leads to \\
\begin{equation}
\left| \sum _{q\geq 1} \frac{\mu (q)}{q^{s-1}\varphi (q)}c_q(n) \right |
\leq 2n \log \log n\sum _{q\geq 1} \frac{\log q}{q^{\sigma }}.
\end{equation}\\

Since for large $q\geq 1$, the inequality \(\log q \leq q^{\varepsilon/2 }\) holds for any small number \(\epsilon >0\), and $\Re e(s)=\sigma \geq1+\varepsilon$,
it is clear that the series is absolutely and uniformly convergent for \(\Re e(s)\geq 1+\varepsilon\). $\blacksquare $\\                                                                      

Detailed analysis of Ramanujan series are given in \cite{HR59}, \cite{KL12}, \cite{LL10}, \cite{RM13}, \cite{SS94}, and similar references.  \\

\subsection{Taylor Series}
The Taylor series $f(n,s)$ at $s=1$ has the ratio $\frac{\varphi(n)}{n}\Lambda(n)=\Lambda(n)+O(\frac{\log n}{n})$ as a constant term. Ergo, this series realizes an effective analytic approximation of the vonMangoldt function $\Lambda(n)$ as a function of two variables $n\geq2$ and $s\approx 1$.\\

\textbf{Lemma 4.3.} \textit{For any integer $n \geq 2$, and any complex number $s \in \mathbb{C}$, $\Re e(s)> 1$, the Taylor series of the function $f(n,s)$ at $s = 1$ is as follows.}\\
\begin{equation}
f(n,s)= f_0(s-1)^{\omega-1}+f_{1}(s-1)^{\omega}+f_{2}(s-1)^{\omega+1}+ \cdot\cdot \cdot,
\end{equation}\\
\textit{where the $kth$ coefficient $f_{k}=f_{k}(n)$ is a function of $n\geq 2$, and $\omega=\omega(n)=\#\{p|n\}$.}\\

\textbf{Proof:} As the series has multiplicative coefficients, it has a product expansion:\\
\begin{eqnarray} \label{el43}
\sum _{q\geq 1} \frac{\mu (q)}{q^{s-1}\varphi
	(q)}c_q(n)
&=&\prod _{p\geq 2} \left(1+\sum _{k\geq 1} \frac{\mu \left(p^k\right)}{p^{k(s-1)}\varphi \left(p^k\right)}c_{p^k}(n)\right)\\
&=&\prod _{p\geq
	2} \left(1+\frac{\mu (p)}{p^{s-1}\varphi (p)}c_p(n)\right)\\
&=&\prod _{p\geq 2} \left(1-\frac{1}{p^{s-1}(p-1)}c_p(n)\right)\\
&=&\prod _{p|n} \left(1-\frac{1}{p^{s-1}}\right)\prod
_{p \nmid n} \left(1+\frac{1}{p^{s-1}(p-1)}\right)\\
&=&\prod _{p|n} \left(\frac{(p-1)\left(p^{s-1}-1\right)}{p^s-p^{s-1}+1}\right)\prod _{p\geq 2} \left(1+\frac{1}{p^{s-1}(p-1)}\right).
\end{eqnarray}\\

These reductions and simplifications have used the identities \(\mu \left(p^k\right)=0,k\geq 2\), in line 2; and \\
\begin{equation}
\begin{array}{ll}  & 
c_p(n)=\left \{\begin{array}{ll}
p-1 & \text{if } \gcd (n,p)=p, \\
-1 & \text{if } \gcd (n,p)=1, \\
\end{array}\right.
\end{array} 
\end{equation}\\
in line 4, see Lemma 2.1. Write the product representation as \\
\begin{equation}
\prod _{p|n} \left(\frac{(p-1)\left(p^{s-1}-1\right)}{p^s-p^{s-1}+1}\right)\prod
_{p\geq 2} \left(1+\frac{1}{p^{s-1}(p-1)}\right) =\zeta (s)A(s)B(s).
\end{equation} \\
The factor \(A(s)=A(n,s)\), a function of two variables of \(n\in \mathbb{N}\) and \(s\in \mathbb{C}\), has the Taylor series \\ 
\begin{eqnarray} \label{el432}
A(s)&=&\prod
_{p|n} \left(\frac{(p-1)\left(p^{s-1}-1\right)}{p^s-p^{s-1}+1}\right) \\ &=&a_1(s-1)^{\omega}+a_2(s-1)^{\omega+1}+a_3(s-1)^{\omega+2}+\cdots,
\end{eqnarray} \\
where the \textit{k}th coefficient \(a_k=a_k(n)\) is a function of \(n\geq 2\), and $\omega=\omega(n)=\#\{p|n\}$ determines the multiplicity of the zero at \(s=1\). \\

\textbf{Case 1:} If \(n=p^v,v\geq 1\), where \(p\geq 2\) is prime, then  $\omega=\omega(n)=1$ and $A(s)$ has a simple zero at \(s=1\). \\

\textbf{Case 2:} If \(n\neq p^v,v\geq 1\), where
\(p\geq 2\) is prime, then $\omega=\omega(n)>1$ and $A(s)$ has a zero of multiplicity $\omega>1$ at \(s=1\). \\

The factor \(B(s)=B(n,s)\), a function of two variables of \(n\in \mathbb{N}\) and \(s\in \mathbb{C}\), has the Taylor series \\
\begin{eqnarray}
B(s)&=&\prod
_{p\geq 2} \left(1+\frac{1}{p^{s-1}(p-1)}\right)\left(1-\frac{1}{p^s}\right)\\
&=&1+b_1(s-1)+b_2(s-1){}^2+b_3(s-1){}^3+\cdots,
\end{eqnarray} \\
which is analytic for \(\Re e(s)>0\). And the zeta function \\
\begin{eqnarray}
\zeta (s)&=&\prod _{p\geq 2} \left(1-\frac{1}{p^s}\right)^{-1}\\
&=&\frac{1}{s-1}+\gamma _0+\gamma _1(s-1)+\gamma _2(s-1)^2+\gamma _3(s-1)^3+\cdots  ,
\end{eqnarray} \\
where \(\gamma _k\) is the \textit{k}th Stieltjes constant, see \cite{KJ92}, and \cite[eq. 25.2.4]{DLMF}. Hence, the Taylor series is \\
\begin{eqnarray} \label{el431}
\sum _{q\geq 1}
\frac{\mu (q)}{q^{s-1}\varphi (q)}c_q(n)
&=&\zeta (s)A(s)B(s) \\
&=&\left(\frac{1}{s-1}+\gamma _0+\gamma _1(s-1)+\gamma _2(s-1)^2+\gamma _3(s-1)^3+\cdots
\right) \\
&\qquad \times& \left(a_1(s-1)^\omega+a_2(s-1)^{\omega+1}+a_3(s-1)^{\omega+2}+\cdots \right) \\
&\qquad \times&\left(1+b_1(s-1)+b_2(s-1)^2+b_3(s-1)^3+\cdots  \right)\\
&=&a_1(s-1)^{\omega-1}+\left(a_1\gamma _0+a_1b_1+a_2 \right)(s-1)^{\omega} \\ 
&\qquad +&\left(a_1b_2+a_1\gamma _1+a_2\gamma _0+a_3+a_1b_1\gamma _0\right)(s-1)^{\omega+1}+\cdots \\
&=&f_0(s-1)^{\omega-1}+f_{1}(s-1)^{\omega}+f_{2}(s-1)^{\omega+1}+ \cdots,
\end{eqnarray}\\
for \(s\in \mathbb{C} \) such that \(\Re e(s) >1\). \\

In synopsis, the Ramanujan series and Taylor series are equivalent: \\
\begin{equation}
\sum _{q\geq 1} \frac{\mu (q)}{q^{s-1}\varphi (q)}c_q(n)=f_0(s-1)^{\omega-1} (n)+f_1(s-1)^{\omega}+f_2(s-1)^{\omega+1}+\cdots .
\end{equation} \\
In particular, evaluating the product representation or the Taylor series or Ramanujan series at \(s=1\) shows that \\
\begin{equation}
 \begin{array}{ll}
& 
f(n,1)=
\left \{\begin{array}{ll}
a_1 & \text{   if  } n=p^k,k\geq 1; \\
0 & \text{   if  } n\neq p^k, k\geq 1; \\
\infty  & \text{   if  } n=1 . \\
\end{array}\right .  \\
\end{array} 
\end{equation} \\
These complete the proof. $\blacksquare$\\

The first few coefficients\\

$f_0=a_1$, \hskip 3 in $f_1=a_1\gamma _0+a_1b_1+a_2$, \\ 
$f_2=a_1b_2+a_1\gamma _1+a_2\gamma _0+a_3+a_1b_1\gamma _0$,\\
 
are determined from the derivatives of \(A(s)\). The first derivative is \\
\begin{equation}
\frac{d}{d s}A(s)=\prod _{p|n} \left(\frac{(p-1)p^{s-1}}{p^s-p^{s-1}+1}-\frac{(p-1)^2\left(p^{s-1}-1\right)p^{s-1}}{\left(p^s-p^{s-1}+1\right)^2}\right)
\log  p .
\end{equation}\\
Accordingly, the first coefficient $a_1=A'(1)$ is written as \\
\begin{equation} \label{el433}
\left. a_1=A^{'}(s) \right|_{s=1}=\prod_{p|n} \left (\frac{p-1}{p}\log p \right ),
\end{equation} 
and the second coefficient $a_2=A^{''}(1)$ is written as \\
\begin{equation} \label{el434}
\left. a_2=A^{''}(s) \right|_{s=1}=\prod_{p|n} \left (\frac{p-1}{p^2}\log^2 p \right ).
\end{equation} \\

\textbf{Example.} This illustrates the effect of the decomposition of the integer $n\geq 1$ on the multiplicity of the zero of the function $f(n,s)$ and some details on the shapes of the coefficients.\\

(i) If $n=p^k,k\geq1$ is a prime power, which is equivalent to $\omega(n)=1$, then\\
\begin{equation}
f(n,s)=f_0+f_1(s-1)+f_2(s-1)^{2}+\cdots,
\end{equation}  
where $f_0=\prod_{p|n} \left (\frac{p-1}{p}\log p \right ) =\frac{\varphi (n)}{n}\Lambda
(n)$ is the most important coefficient. The other coefficients are absorbed in the error term.\\

(ii) If $n=p_1^a p_2^b,a,b\geq 1$ is a product of two prime powers, which is equivalent to $\omega(n)=2$, then\\
\begin{equation}
f(n,s)=f_0(s-1)+f_1(s-1)^{2}+f_2(s-1)^{3}+\cdots,
\end{equation}
where $f_0=\prod_{p|n} \left (\frac{p-1}{p}\log p \right ) \ne\frac{\varphi (n)}{n}\Lambda
(n)$. All the coefficients, including $f_0$, are absorbed in the error term.\\

These information are used in the next result to assemble an effective approximation of the vonMangoldt function:\\  
\begin{equation}
\begin{array}{ll}
& 
f(n,s)=
\left \{\begin{array}{ll}
\frac{ \varphi(n)}{n} \Lambda(n)+O(\varepsilon) & \text{   if  } n=p^k,k\geq 1, \text{ and } 0<|s-1|<\delta; \\
O(\varepsilon) & \text{   if  } n\neq p^k, k\geq 1, \text{ and } 0<|s-1|<\delta. \\
\end{array}\right .  \\
\end{array} 
\end{equation} 
The details and the proof are given below.\\
                                 
\textbf{Lemma 4.4.}   \textit{Let $\varepsilon >0$ be an arbitrarily small number, and let $n \geq 2$. Then, the Taylor series satisfies} \\
\begin{enumerate}[(i)]
\item $\displaystyle \left | f(n,s)-\frac{ \varphi(n)}{n} \Lambda(n) \right | \leq \varepsilon, $ 
\item $\displaystyle f(n,s)=\frac{ \varphi(n)}{n} \Lambda(n) +O(\varepsilon), $
\item $\displaystyle f(n,s)=\frac{ \varphi(n)}{n} \Lambda(n) +f^{'}(n,s_0)(s-1), $ \textit{for some real number} $0<|s_0-1|<\delta$, 
\end{enumerate}	
\textit{where $0< |s-1|< \delta$, and a constant $\delta=c \varepsilon>0$, with $c>0$ constant.}\\

\textbf{Proof:} For $n > 1$, the Taylor series $f(n,s)$ is a continuous function of a complex number $s \in \mathbb{C}, \Re e(s) = \sigma>1$. It has continuous and absolutely convergent derivatives\\
\begin{equation}
f^{(k)}(n,s)=(-1)^k \sum _{q\geq 1} \frac{\mu (q)\log^k q}{q^{s-1}\varphi (q)}c_q(n),
\end{equation}
 on the half plane $\Re e(s)>1$ for $k \geq1$ . By the definition of continuity, for any $\varepsilon>0$ arbitrarily small number, there exists a small real number $\delta>0$ such that\\
\begin{equation}
\left | f(n,s)-f(n,1) \right | < \varepsilon
\end{equation}
for all $s \in \mathbb{C}$ such that $0< |s-1|< \delta$. The relation $\delta=c \varepsilon>0$, with $c>0$ constant, is a direct consequence of uniform continuity, confer Lemma 4.2. \\

The statement (iii) follows from the Mean Value Theorem:\\
\begin{equation}
 (s-1) f^{'}(n,s_0) =f(n,s)-f(n,1)
\end{equation}\\
for some fixed real number $0<|s_0-1|<\delta$. Moreover, choice $\delta=\varepsilon/(2f^{'}(n,s_0))$ is sufficient. $\blacksquare$ \\

The uniform convergence of the derivative $f^{'}(n,s)$ for $\Re e(s)>1$, and the Mean Value Theorem imply that the product $|(s-1)f^{'}(n,s)|< \varepsilon $ is uniformly bounded for all $n\geq 1$ and $\Re e(s)>1$. Accordingly, statements (ii) and (iii) in Lemma 4.4 are equivalent, and can be interchanged in the proof of Lemma 6.1.\\

There is some interest in obtaining additional information on of the properties of the derivative, (not required in the proof of Theorem 1.1). The next result shows that the derivative $f^{'}(n,s)$ is zerofree on the half plane $\{s \in \mathbb{C}: \Re e(s)>1\}$ independently of $n \geq 1$. \\

\textbf{Lemma 4.5.} \textit{For any integer $n \geq 1$, and any complex number $s \in \mathbb{C}$, $\Re e(s)> 1$, the derivative of 
the function $f(n,s)$ satisfies $ f^{'}(n,s) >0$ for $\Re e(s)>1$}.\\

\textbf{Proof:} Start with the fourth line in equation (\ref{el43}), that is,\\
\begin{equation}
f(n,s)=\sum _{q\geq 1} \frac{\mu (q)}{q^{s-1}\varphi(q)}c_q(n)=\prod _{p|n} \left(1-\frac{1}{p^{s-1}}\right) \prod
_{p \nmid n} \left(1+\frac{1}{p^{s-1}(p-1)}\right).
\end{equation}\\ 
The corresponding derivative is\\ 
\begin{eqnarray}
f^{'}(n,s)&=&\prod _{p|n} \left(1+\frac{\log p}{p^{s-1}}\right)\prod
_{p \nmid n} \left(1+\frac{1}{p^{s-1}(p-1)}\right)\\
&+& \prod _{p|n} \left(1-\frac{1}{p^{s-1}}\right)\prod
_{p \nmid n} \left(1-\frac{\log p}{p^{s-1}(p-1)}\right).
\end{eqnarray}\\ 
Now, observe that the first double product is positive and larger than the second double product. Hence $f^{'}(n,s)>0$ for real $s>1$.    $\blacksquare$\\

In summary, both Ramanujan series and the Taylor series for $f(n,s)$ are effective approximations of the vonMangoldt function, id est, \\
\begin{equation}
 f(n,s)=\frac{ \varphi(n)}{n} \Lambda(n) +f^{'}(n,s_0)(s-1) \approx \Lambda(n), 
\end{equation}
for $n\geq2$ if $|s-1|<\delta$, with $\delta>0$ a small number, and for some fixed real number $0<|s_0-1|<\delta$. This demonstrated and proved in Lemma 4.4-iii.\\

Many other related approximations and identities for the vonMangoldt function appear in \cite[pp. 25--28]{HG07}, \cite{GG05}, \cite{FI10}, \cite{KL12}, and the literature.\\

\vskip 1 in
\section{The Wiener Khintchine Formula}
A basic result in Fourier analysis on the correlation function of a pair of continuous functions, known as the Wiener-Khintchine formula, is an important part of this work.  
\\

\textbf{Theorem 5.1.}  \textit{If $f(s)=\sum_{n\geq1} a_{q} c_{q}(n) $ is an absolutely convergent Ramanujan series of a number theoretical function $f : \mathbb{N}\longrightarrow \mathbb{C}$, and $m \geq 1$ is a fixed integer, then, }
\\
\begin{eqnarray}
\lim_{x\rightarrow \infty} \frac{1}{x} \sum_{n\leq x}f(n,s)f(n+m,s)=\sum_{q\geq 1} a_{q}^{2}c_{q}(m).
\end{eqnarray}
\\
\textbf{Proof:} Compute the mean value of the product $f(n,s)f(n+m,s)$ as stated in Section 3. This requires the Ramanujan series, and the exponential sum $c_q(n)=\sum _{\gcd(u,q)=1} e^{i2\pi nu/q} $ as demonstrated here.\\
\begin{eqnarray*}
\lim_{x\longrightarrow \infty } \frac{1}{x}\sum _{n\leq x} f(n)f(n+m)
&=&\lim_{x\longrightarrow \infty } \frac{1}{x}\sum _{n\leq x} \left(\sum _{q\geq
	1} a_qc_q(n)\right)\left(\sum _{r\geq 1} a_rc_r(n+m)\right)\\
&=&\lim_{x\longrightarrow \infty } \frac{1}{x}\sum _{n\leq x} \left(\sum _{q\geq 1} a_q\sum _{\gcd(u,q)=1} e^{\text{i2$\pi $} n \frac{u
	}{q}}\right)\left(\sum _{r\geq 1} a_r\sum _{\gcd (v,r)=1} e^{\text{i2$\pi $} (n+m) \frac{v }{r}}\right).
\end{eqnarray*}
\begin{equation}
\quad
\end{equation}\\
The last equation follows from the orthogonal mean value relation, Lemma 3.1. Next, use the absolute convergence to exchange the order of summation:
\begin{eqnarray*}
\lim_{x\rightarrow \infty } \frac{1}{x}\sum _{n\leq x} f(n)f(n+m)
&=&\sum _{q\geq 1} \sum _{r\geq 1} a_qa_r\lim_{x\rightarrow \infty
} \frac{1}{x}\sum _{n\leq x} \sum _{\gcd (u,q)=1} e^{i2\pi n \frac{u }{q}}\sum_{\gcd (v,r)=1} e^{-i2\pi 
(n+m) \frac{v}{r}}\\
&=&\sum_{q\geq 1} \sum_{r\geq 1} a_qa_r \sum _{\gcd(v,r)=1} e^{i2\pi m\frac{v }{r}}\lim_{x\rightarrow \infty}\sum _{n\leq x} \frac{1}{x}\sum_{\gcd (u,q)=1} e^{-i2\pi
n\left (\frac{u}{q}- \frac{v }{r}\right)}\\
&=&
\sum _{q\geq 1} a_qa_qc_q(m).
\end{eqnarray*}
\begin{equation}
 \quad   
\end{equation}\\
Reapply Lemma 3.1 to complete the proof.  $\blacksquare$\\

The Wiener-Khintchine Theorem is related to the Convolution Theorem in Fourier analysis. Other related topics are Parseval formula, Poisson summation formula, et cetera, confer the literature for more details. Some earliest works are given in \cite{WN66}.
\\

A new result on the error term associated with the correlation of a pair of functions complements the well established Wiener-Khintchine formula. This is an essential part for applications in number theory.\\

\textbf{Theorem 5.2.}  (\cite{CM15}) \textit{Let $f$ and $g$ be two arithmetic functions with absolutely convergent Ramanujan series}\\
\begin{equation}
f(s)=\sum_{q\geq1} a_{q} c_{q}(n) \;    \mathrm{and} \;    g(s)=\sum_{q\geq1} b_{q} c_{q}(n),
\end{equation}\\
\textit{where the coefficients satisfy $a_{q},b_{q}=O(q^{-1-\beta})$, with $0 < \beta < 1$, and let $m\geq1$ be a fixed integer. Then,} \\
\begin{equation}
\sum_{n\leq x} f(n,s)f(n+m,s) 
=x \sum_{q\geq1}a_{q}b_{q}c_{q}(m)+O\left(x^{1-\beta}(\log x)^{4-2\beta}\right) .
\end{equation}\\

For $\beta>1$, it readily follows that the error term is significantly better. The proof of this result is lengthy, and should studied in the original source.\\

\vskip 1 in
\section{Approximations For The Correlation Function}
Analytic formulae for approximating the autocorrelation function $\sum_{n \leq x} f(n,s) f(n+m,s)$ are considered here. Two different ways of approximating the correlation function of $f(n,s)$ are demonstrated.\\

\subsection{Taylor Domain}
The Taylor series is used to derive an analytic formula for approximating the correlation function of $f(n,s)$ in terms of the function $\Lambda(n)$ and its shift $\Lambda(n+m)$.\\

\textbf{Lemma 6.1.} \textit{Let $x\geq 1$ be a large number, and let $m\geq 1$ be a fixed integer. Then, the correlation function satisfies}\\
\begin{equation}
\sum_{n\leq x} f(n,s)f(n+m,s) 
= \sum_{n\leq x}\Lambda(n)\Lambda(n+m)+O \left (x\frac{(\log\log x)^2}{\log x} \right ),
\end{equation}\\
\textit{where $0<|s-1|<\delta$, and $0< \delta \ll ({\log \log x)^2}/{\log x}$.}\\

\textbf{Proof:} Let $x \geq 1$ be a large number, and let  $\varepsilon=(\log \log x)^2/\log x>0$. By Lemma 4.4, the correlation function can be rewritten as\\ 
\begin{eqnarray}  \label{el40}
\sum_{n\leq x} f(n,s)f(n+m,s) 
&=& \sum_{n\leq x}\left (\frac{\varphi (n)}{n}\Lambda
(n)+O(\varepsilon) \right )\left (\frac{\varphi (n+m)}{n+m}\Lambda(n+m)+O(\varepsilon) \right )\\
&=&\sum_{n\leq x}\frac{\varphi (n)}{n}\Lambda(n)\frac{\varphi (n+m)}{n+m}\Lambda
(n+m) \\
& &+O\left(\varepsilon\sum_{n\leq x} \left (\frac{\varphi (n)}{n}\Lambda(n)+\frac{\varphi (n+m)}{n+m}\Lambda
(n+m) +\varepsilon\right ) \right )\\
&=&\text{Main Term+Error Term},
\end{eqnarray}\\
where $0<|s-1|<\delta$ for some small number $0< \delta \ll (\log \log x)^2/\log x=\varepsilon$. Now, replace\\
\begin{equation}
\frac{\varphi (n)}{n}\Lambda
(n)=\Lambda
(n)+O \left (\frac{\log n}{n}\right),
\end{equation}\\
and \\
\begin{equation}
\sum_{n\leq x}\frac{\Lambda(n)}{n}=\log x +O(1),
\end{equation}\\
see \cite[p. 88]{AP76}, and \cite[p. 348]{HW08}, into the main term to reduce it to\\
\begin{eqnarray*}
\sum_{n \leq x} \frac{\varphi (n)}{n} \Lambda(n) \frac{\varphi (n+m)}{n+m} \Lambda(n+m)
&=& \sum_{n \leq x}\left (\Lambda(n)+O\left (\frac{\log n}{n}\right) \right ) \left(\Lambda(n+m)+O\left (\frac{\log n}{n} \right ) \right )\\
&=& \sum_{n \leq x}\Lambda(n)\Lambda(n+m)\\
&  &+O\left (\sum_{n \leq x}\left (\frac{\Lambda(n)\log n}{n}+\frac{\Lambda(n+m)\log n}{n}+ \frac{\log^2 n}{n^2}\right ) \right ) \\
&= & \sum_{n \leq x} \Lambda(n) \Lambda(n+m)+O\left ( \log^2 x \right ).
\end{eqnarray*}
\begin{equation} \label{el41}
\quad       
\end{equation}
Note that for fixed $m\geq 1$, and large $n\geq 2$, the expression
\begin{equation} 
\frac{\varphi(n+m)}{n+m}\Lambda(n+m)=\Lambda(n+m)+
O(\frac{\log n}{n} )       
\end{equation}\\
was used in the previous calculations. Similarly, apply the Prime Number Theorem\\
\begin{equation}
\sum_{n\leq x}\Lambda(n)=x +O\left (x e^{-c \sqrt{\log x}} \right),
\end{equation}\\
where $c>0$ is an absolute constant, refer to \cite[p. 179]{MV07}, \cite[p. 277]{TG15}, into the error term to reduce it to\\
\begin{eqnarray} \label{el42}
\sum_{n\leq x} \left (\frac{\varphi (n)}{n}\Lambda(n)+\frac{\varphi (n+m)}{n+m}\Lambda(n+m)+\varepsilon \right )
&=&\sum_{n\leq x}\left (\Lambda(n)+O\left(\frac{\log n}{n}\right) \right ) \\
&  & +\sum_{n\leq x}\left (\Lambda(n+m)+O\left(\frac{\log n}{n}\right)\right )+\varepsilon x\\
&=&2x+\varepsilon x+O\left (x e^{-c \sqrt{\log x}} \right) \\&=&O(x).
\end{eqnarray}\\
Replace (\ref{el41}) and (\ref{el42}) back into (\ref{el40}) to obtain\\
\begin{eqnarray}  \label{el44}
\sum_{n\leq x} f(n,s)f(n+m,s) 
&=& \text{Main Term+Error Term}\\
&=& \sum_{n \leq x} \Lambda(n) \Lambda(n+m)+O\left ( \log^2 x \right )+O(\varepsilon x)\\
&=& \sum_{n \leq x} \Lambda(n) \Lambda(n+m)+O \left (x\frac{(\log\log x)^2}{\log x} \right ).
\end{eqnarray}\\
These confirm the claim.     $ \blacksquare$ \\

\subsection{Ramanujan Domain}
The Ramanujan series is used to derive an analytic formula for approximating the correlation function of $f(n,s)$ in terms of the constant $\alpha_m=\sum_{q\geq1}\left | \frac{\mu(q)}{q^{s-1}\varphi(q)} \right |^{2}c_{q}(m) \geq 0$.\\

\textbf{Lemma 6.2.} \textit{Let $x\geq 2$ be a large number, and let $m\geq 1$ be a fixed integer. Then, the correlation function satisfies}\\
\begin{equation}
\sum_{n\leq x} f(n,s)f(n+m,s) 
=\alpha_{m}x+O \left(x\frac{(\log\log x)^{2})}{\log x}\right),
\end{equation}\\
\textit{where $0<|s-1|<\delta$, with $0< \delta \ll ({\log \log x)^2}/{\log x}$.}\\

\textbf{Proof:} Applying Theorem 5.1,  (Weiner-Khintchine formula), yields the spectrum (mean value of the correlation function of $f(n,s)$): \\
\begin{equation}
\lim_{x\rightarrow \infty} \frac{1}{x} \sum_{n \leq x}f(n,s)f(n+m,s)=\sum_{q\geq 1} \left | \frac{\mu(q)}{q^{s-1}\varphi(q)} \right |^{2}c_{q}(m)=\alpha_{m}.
\end{equation}\\
For fixed $m \geq 1$, the nonnegative constant $\alpha_{m} \geq 0$ determines the density of the subset of prime pairs $\{ p, p + m: p\text{ is prime} \} \subset \mathbb{P}$ in the set of primes $\mathbb{P} = \{ 2, 3, 5,\ldots \}$, see Theorem 7.1 for more details.\\

By Lemma 4.2, the \textit{q}th coefficient satisfies
\\
\begin{equation}
\left | a_{q}(s)\right | = \left | \frac{\mu(q)}{q^{s-1}\varphi(q)} \right | \leq\frac{\log q}{q^{\sigma}}\leq\frac{1}{q^{1+\beta}},
\end{equation} \\
where $\beta=c\varepsilon/2=c({\log \log x)^2}/{ 2 \log x}>0$.\\

Applying Theorem 5.2, to the correlation function of the Ramanujan series of $f(n,s)$ returns the asymptotic formula\\
\begin{eqnarray} 
\sum_{n\leq x} f(n,s)f(n+m,s) 
&=&x \sum_{q\geq1}\left | \frac{\mu(q)}{q^{s-1}\varphi(q)} \right |^{2}c_{q}(m)+O(x^{1-\beta}(\log x)^{4-2\beta})\\
&=& \alpha_{m} x+O \left (x\frac{(\log\log x)^2}{\log x} \right ),
\end{eqnarray}
\\
where $\Re e(s)=1+c(\log\log x)^{2} / \log x=1+2\beta>1$ as required by Theorem 5.2, and 
\\ 
\begin{equation}
x^{1-\beta}(\log x)^{4-2\beta}=x^{1-\frac{(\log\log x)^2}{ 2\log x}} (\log x)^{4} \leq x \frac{(\log\log x)^2}{\log x}.
\end{equation}\\
This holds for all complex numbers $s \in \mathbb{C}$ for which $0<|s-1|<\delta$, with $0< \delta \ll \varepsilon\leq (\log \log x)^{2}/{\log x}$, see Lemma 4.4.  $\blacksquare$\\

\vskip 1 in
\section{Spectrum of $f(n,s)$}
For a fixed integer $m\geq1$, the constant $\alpha_m$ approximates the density of the subset of prime pairs $\{p,p+m\}\subset \mathbb{P}$ in the set of primes $\mathbb{P}=\{2,3,5,7, \ldots\}$. Some of the basic properties of the densities associated with the subset of primes pairs are explained in this Section.\\

\subsection{The Function $F(\sigma)$}
The coefficients of the series $f(n,s)$ are multiplicative. This property enable simple derivation and simple expressions for the constants, which approximate or determine the densities of the subsets for the various primes counting problems. Furthermore, the spectrum analysis in Theorem 7.1 provides a detailed view of certain properties of the densities  as a function of $m \geq 1$ and $\sigma>1/2$.\\

\textbf{Theorem 7.1.}  Let $m \geq1$ be a fixed integer, and let $\Re e(s)=\sigma>1$ be a complex number. Then, the spectrum of the correlation function $f(n,s)$ has a product decomposition as\\
\begin{equation}
F(\sigma)=\prod_{p|m} \left (1+\frac{1}{p^{2\sigma-2}(p-1)}\right ) \prod_{p\nmid m} \left (1-\frac{1}{p^{2\sigma-2}(p-1)^{2}}\right ) .
\end{equation} \\
In particular, the followings statements hold.\\

(i) If $m=2k+1$ is a fixed odd integer, then the density of prime pairs is $F(1)=0$. This implies that the number of prime pairs is finite or a subset of zero density.\\

(ii) If $m=2k\geq 2$ is a fixed even integer, and $\Re e(s)=\sigma \geq 1$, then the density of prime pairs is $F(\sigma)\geq 1/2$. This implies that the number of prime pairs is infinite.\\

\textbf{Proof:} By Lemma 4.2, the function $f(n,s)$ has an absolutely and uniformly convergent Ramanujan series for $\Re e(s)=\sigma>1$. Applying Theorem 5.1 (Weiner-Khintchine formula) yields the spectrum\\
\begin{eqnarray}
F(\sigma)&=& \lim_{x \rightarrow \infty}\frac{1}{x} \sum_{n \leq x} f(n,s)f(n+m,\overline{s})\\
&=&\lim_{x \rightarrow \infty}\frac{1}{x} \sum_{n \leq x} \left ( \sum_{q\geq 1} \frac{\mu(q)}{q^{s-1}\varphi(q)} c_q(n)\right)\left ( \sum_{r\geq 1} \frac{\mu(r)}{r^{s-1}\varphi(r)} c_r(n+m)\right)    \\
&=&\sum_{q\geq 1} \left | \frac{\mu(q)}{q^{s-1}\varphi(q)} \right |^{2} c_{q}(m).
\end{eqnarray}\\

Using the multiplicative property of the coefficients returns the conversion to a product\\
\begin{equation}
F(\sigma)=\sum_{q\geq 1} \left | \frac{\mu(q)}{q^{s-1}\varphi(q)} \right |^{2} c_{q}(m)=\prod_{p\geq 2} \left (1+\frac{\mu(p)^2}{p^{2\sigma-2}\varphi(p)^2} \overline{c}_{p}(m)\right ) .
\end{equation}\\
If $m=2k+1$ is odd, then the factor for $q=2$ is\\
\begin{equation}
1+\frac{\mu(p)^2}{p^{2\sigma-2}\varphi(p)^2} \overline{c}_{p}(m)=1-\frac{1}{2^{2\sigma-2}}.
\end{equation}\\
This follows from Lemma 2.1, that is, $c_q(m)=\overline{c_p(m)}$ and \\
\begin{equation}
c_{p}(m)=\left \{\begin{array}{ll}
p-1 & \text{if } \gcd (m,p)=p, \\
-1 & \text{if } \gcd (m,p)=1. 
\end{array}\right.
\end{equation} \\
Consequently, the spectrum $F(\sigma)$ vanishes at $s=\sigma=1$, this proves (i). While for even $m=2k\geq 2$, it reduces to\\
\begin{equation}
F(\sigma)=\prod_{p |2k} \left (1+\frac{1}{p^{2\sigma-2}(p-1)}\right ) \prod_{p \not 2k} \left (1-\frac{1}{p^{2\sigma-2}(p-1)^{2}}\right ) .
\end{equation}\\
Lastly, it clear that the product converges and $F(\sigma)>0$ for $\Re e(s)=\sigma>1$.  $\blacksquare$\\

\subsection{Numerical Data For $F(\sigma)$}
Some numerical data was compiled for the spectrum function $F(\sigma)$ of the correlation function of $f(n,s)$. The computation uses $F(\sigma)= \left (1+\frac{1}{2^{2\sigma-2}}\right ) \prod_{3\leq p \leq x} \left (1-\frac{1}{p^{2\sigma-2}(p-1)^{2}}\right ) $, and it was restricted to the primes in the interval $p\leq x=10^4 \log 10^4$ and $1/4 \leq \sigma \leq 10$.\\

\begin{figure}[h]
\centering
\includegraphics[width=24cm,height=8cm]{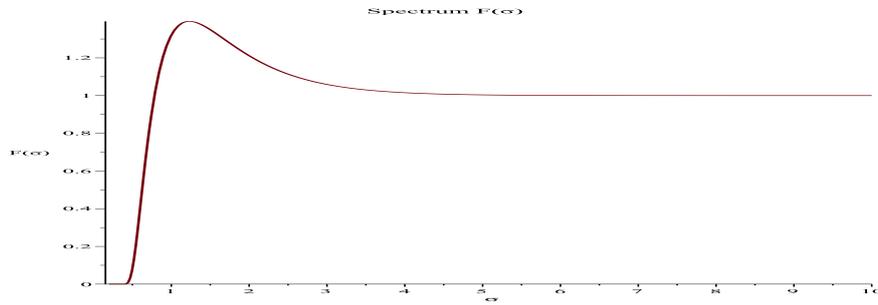}
\caption{Graph of the Density of Twin Primes}
\label{fig:mesh1}
\end{figure}

The value $F(1)=1.3205814800...$ corresponds to the twin prime constant. The peak $F(1.2359486...)=1.3894427...$ appears to be unexpected. This peak probably heralds certain irregularity in the distribution of twin primes. 

\subsection{The Limits $F(1)$}
The constant $\alpha_k>0$ in Theorem 1.1 is given by a fixed value of the spectrum function, that is, $\alpha_k=F(\sigma),\sigma>1$ fixed. For $1<\sigma$, it is clear that the product converges, and $\alpha_k>1/2$. This demonstrates that the density of prime pairs $p,p+2k$ increases as the parameter $\sigma\rightarrow 1^{+}$. In particular, as  $ \sigma=1+2\beta=1+2c\varepsilon>1$ decreases, the limit of the spectrum is\\
\begin{eqnarray}
A_{k}&=&\lim_{\varepsilon\rightarrow 0} F(1+c\varepsilon)\\
&=&\lim_{\varepsilon\rightarrow 0}\prod_{p| 2k} \left (1+\frac{1}{p^{2c\varepsilon}(p-1)}\right ) \prod_{p \nmid 2k} \left (1-\frac{1}{p^{2c\varepsilon}(p-1)^{2}}\right )\\
&=& \prod_{p| 2k} \left (1+\frac{1}{p-1}\right ) \prod_{p\nmid 2k} \left (1-\frac{1}{(p-1)^{2}}\right ),
\end{eqnarray}\\
refer to \cite{RP96}, \cite{NW00}, \cite{GG05}, et alii, for more information. This is precisely the constant\\
\begin{equation}
A_{k}=\prod_{p|2k} \left (1+\frac{1}{p-1}\right ) \prod_{p\nmid 2k} \left (1-\frac{1}{(p-1)^{2}}\right )\geq \alpha_{k}. 
\end{equation}
associated with the density of the subset of prime pairs $\{p,p+2k\}$ in the set of primes $\mathbb{P} = \{ 2, 3, 5, 7, 11, \ldots \}$. For example, for $2k = 2$, this is the twin prime constant:\\
\begin{equation}
A_{2}=2 \prod_{p\geq3} \left (1-\frac{1}{(p-1)^2}\right )=1.32058148001344\ldots \:.
\end{equation}\\

\subsection{Circle Methods}
The circle method is the standard analytic tool used to derive these constants. The constants (or equivalently the densities) for several primes counting problems are given in terms of some singular series such as\\
\begin{equation}
\mathfrak{G}_{r}(m)=\sum_{q \geq 1} \left (\frac{\mu(q)}{\varphi(q)}\right )^{r}c_q(m), 
\end{equation}\\
where $m,r \geq 1$. Some of the earliest works appears in \cite{HL23}. More recent works, at various levels of difficulties, are given in \cite{GP06}, \cite[p. 295]{EL85}, \cite[Chapter 23]{MT06}, and \cite{VR97}. \\ 

The specific case for the number of representation $r(N)=(1/2)\mathfrak{G}(N) x+O(N^2 \log ^{-B}N)$ of a large odd integer as a sum $N=p_1+p_2+p_3$ of three primes has the constant\\
\begin{equation}
\mathfrak{G}(N)=\sum_{q \geq 1} \frac{\mu(q)}{\varphi(q)^3}c_q(N)=\prod_{p|N} \left (1-\frac{1}{(p-1)^2}\right ) \prod_{p\nmid N} \left (1+\frac{1}{(p-1)^{3}}\right )>0. 
\end{equation}
The proof appears in \cite[Chapter 26]{DH00}.\\

\subsection{Dirichlet Density}
The Dirichlet density for the subset of primes $\mathcal{P}(a,q)=\#\{p=qn+a:n \geq 1\} \subset \mathbb{P}$ in an arithmetic progression $\{qn+a:\gcd(a,q)=1 \text{ and } n\geq 1\}$ is also defined by a limit. This has the form\\
\begin{equation}
\delta(\mathcal{P}(a,q))=\lim_{s\rightarrow 1^+}\frac{\sum_{p \equiv a \text{ mod } q}p^{-s}} {\sum_{p \geq 2}p^{-s}}=\dfrac{1}{\varphi(q)}.
\end{equation}\\
A proof appears in \cite[p. 237]{CO07}. \\

\vskip 1 in
\section{Linear Twin Primes And Prime Pairs}
The qualitative version of the prime pairs is attributed to dePolinac, \cite{PJ09}, and the heuristic for the quantitative prime pairs conjecture, based on the circle method, appears in \cite[p. 42]{HL23}. The precise statement is the following.\\ 

\textbf{Prime Pairs Conjecture.}  There are infinitely many prime pairs $p, p+2k$ as the prime $p \to \infty$. Moreover, the counting function has the asymptotic formula\\
\begin{equation}
\pi _{k}(x)=2 \prod_{p |k} \left( \frac{p-1}{p-2} \right ) \prod _{p\geq 3} \left(1-\frac{1}{(p-1)^2}\right) \frac{x}{\log ^2 x}+O\left(\frac{x}{\log ^3 x}\right),
\end{equation}\\
where \(k \geq 1\) is a fixed integer, and $x \geq 1$ is a large number.\\

The heuristic based on probability, yields the Gaussian type analytic formula\\
\begin{equation}
\pi _{k}(x)=C_k \int_{2}^{x} \frac{1}{\log^2 t} dt+O\left(\frac{x}{\log ^3 x}\right),
\end{equation}\\
where $C_k$ is the same constant as above. The deterministic approach is based on the weighted prime pairs counting function\\
\begin{equation}
\sum_{n\leq x}\Lambda(n)\Lambda(n+m),
\end{equation}\\
which is basically an extended version of the Chebishev method for counting primes.\\

For an odd integer $m\geq 1$, the average order of the finite sum $\sum_{n\leq x}\Lambda(n)\Lambda(n+m)$ is very small, and the number of prime pairs $p, p + m$ is finite. This is demonstrated in Theorem 7.1. Another way of estimating this is shown in \cite[p. 507]{GP99}.  Consequently, it is sufficient to consider even integer $m = 2k$.
\\

\subsection{The Proof Of Theorem 1.1}
The proof of Theorem 1.1 exploits the duality or equivalence of the Ramanujan series domain and the Taylor series domain to compute the correlation function of $f(n,s)$ in two distinct ways. \\

\textbf{Proof of Theorem 1.1:} Let $x \geq 1$ be a large number, and let $\varepsilon= ({\log \log x)^2}/{\log x}>0$. By Lemma 4.2, the function $f(n,s)$ has an absolutely, and uniformly convergent Ramanujan series for 
$\Re e(s) = \sigma \geq 1+c({\log \log x)^2}/{\log x}>1$, with $0<c<1$ constant. 
\\

By Lemma 6.1, the correlation function of the Taylor series of $f(n,s)$ has the asymptotic formula
\\
\begin{equation} \label{el70} 
\sum_{n\leq x} f(n,s)f(n+2k,s) 
= \sum_{n\leq x}\Lambda(n)\Lambda(n+2k)+O \left (x\frac{(\log\log x)^2}{\log x} \right ),
\end{equation}
\\
where $0<|s-1|<\delta$, with $0< \delta \ll \varepsilon\leq (\log \log x)^{2}/{\log x}$. \\

By Lemma 6.2, the correlation function of the Ramanujan series of $f(n,s)$ has the asymptotic formula\\
\begin{equation} \label{el71}
\sum_{n \leq x}f(n,s)f(n+2k,s)=\alpha_{k}x+O \left(x\frac{(\log\log x)^{2}}{\log x}\right),
\end{equation}\\
Since the constant $\alpha_k>1/2$ for $\Re e(s)>1$, see Theorem 7.1, combining (\ref{el70}) and (\ref{el71}) lead to the inequalities\\
\begin{eqnarray}
\sum_{n\leq x}\Lambda(n)\Lambda(n+2k)+O \left (x\frac{(\log\log x)^2}{\log x} \right )
&=& \sum_{n \leq x}f(n,s)f(n+2k,s)\\
&=& \alpha_{k}x+O \left(x\frac{(\log\log x)^{2}}{\log x}\right)\\
&\geq &  (1/2)x+O \left(x\frac{(\log\log x)^{2}}{\log x}\right).
\end{eqnarray}
These confirm the claim. $\blacksquare$\\

The nonnegative constant\\
\begin{equation}
\alpha_2=\sum_{q\geq 1} \left | \frac{\mu(q)}{q^{s-1}\varphi(q)} \right |^2c_{q}(2)\geq \frac{1}{2}
\end{equation} \\
for $\sigma>1$ shows that the subset of prime pairs $\{ p, p+2k: p\geq 2 \text{ prime}\}$ has positive density in the set of primes $\mathbb{P} = \{ 2,3,5,7,11,\ldots\}$.
\\

\vskip 1 in
\section{Quadratic Twin Primes}
The quadratic twin primes conjecture calls for an infinite subset of prime pairs $n^2+1, n^2+3$ as the integer $n \geq 1$ tends to infinity. Background information, and various descriptions of this primes counting problem are given in \cite[p. 46]{HL23}, \cite[p. 343]{NW00}, \cite[p. 406]{RP96},  \cite{PJ09}, et alii.
\\

The original derivation of the heuristic, based on the circle method, which appears in \cite[p. 62]{HL23}, concluded in the followings observation.\\ 

\textbf{Quadratic Twin Primes Conjecture.}  There are infinitely many twin primes $n^2+1, n^2+3$ as $n \to \infty$. Moreover, the counting function has the asymptotic formula\\
\begin{equation}
\pi _{f_2}(x)=3\frac{x^{1/2}}{\log ^2 x}\prod _{p\geq 5} \left(\frac{p(p-\varrho (p))}{(p-1)^2}\right)+O\left(\frac{x^{1/2}}{\log ^3 x}\right),
\end{equation}\\
where \(\varrho (p)=0,2,4\) according to the occurrence of $0,1,\text{ or }2$ quadratic residues $\text{mod } p$ in the set $\{-1,-3\}$.\\

The heuristic based on probability, yields the Gaussian type analytic formula\\
\begin{equation}
\pi _{f_2}(x)=s_{2} \int_{2}^{x} \frac{1}{ \sqrt t \log^2 t} dt+O\left(\frac{x^{1/2}}{\log ^3 x}\right),
\end{equation}\\
where $s_{2}$ is the same constant as above. The deterministic approach is based on the weighted prime pairs counting function\\
\begin{equation}
\sum_{n\leq x}\Lambda(n^2+1)\Lambda(n^2+3),
\end{equation}\\
which is basically an extended version of the Chebishev method for counting primes.\\

\subsection{The Proof of Theorem 1.2}
A lower bound for the correlation function $\sum_{n^{2}+1\leq x}\Lambda(n^{2}+1)\Lambda(n^{2}+3)$ will be derived here. It exploits the duality or equivalence of the Ramanujan series domain and the Taylor series domain to compute the correlation function of in two distinct ways. \\

\textbf{Proof of Theorem 1.2:} Let $x \geq 1$ be a large number, and let $\varepsilon= ({\log \log x)^2}/{\log x}>0$. By Lemma 2.1, the function $f(n^{2}+1,s)$ has an absolutely, and uniformly convergent Ramanujan series for 
$\Re e(s) = \sigma \geq 1+c({\log \log x)^2}/{\log x}>1$, with $0<c<1$ constant. \\

By Lemma 6.1, the correlation function of the Taylor series of $f(n^{2}+1,s)$ has the asymptotic formula\\
\begin{equation} \label{el90}
\sum_{n^{2}+1\leq x} f(n^{2}+1,s)f(n^{2}+3,s) 
= \sum_{n^{2}+1\leq x}\Lambda(n^{2}+1)\Lambda(n^{2}+3)+O \left (x^{1/2}\frac{(\log\log x)^2}{\log x} \right ),
\end{equation}\\
where $0<|s-1|<\delta$, with $0< \delta \ll \varepsilon=({\log \log x)^2}/{\log x}$, see Lemma 4.4.\\

Applying Lemma 6.2, to the correlation function of the Ramanujan series of $f(n^{2}+1,s)$ returns the asymptotic formula\\
\begin{equation} \label{el91}
\sum_{n^{2}+1\leq x} f(n^{2}+1,s)f(n^{2}+3,s) 
=\alpha_{2} x^{1/2}+O \left (x^{1/2}\frac{(\log\log x)^2}{\log x} \right ),
\end{equation}
\\
where $\Re e(s)=1+c(\log\log x)^{2}/\log x=1+2\beta>1$ as required by Theorem 6.2. \\ 

Since the constant $\alpha_k>1/2$ for $\Re e(s)>1$, see Theorem 7.1, combining (\ref{el90}) and (\ref{el91}) confirms the claim: \\
\begin{eqnarray}
\sum_{n^2+1\leq x}\Lambda(n^{2}+1)\Lambda(n^{2}+3) 
&=& \sum_{n^{2}+1\leq x} f(n^{2}+1,s)f(n^{2}+3,s)  \\
&=& \alpha_{2}x^{1/2}+O \left(x^{1/2}\frac{(\log\log x)^{2}}{\log x}\right)\\
&\geq& (1/2)x^{1/2}+O \left(x^{1/2}\frac{(\log\log x)^{2}}{\log x}\right),
\end{eqnarray}\\
for a lower bound of the correlation function.   $\blacksquare$  \\

The nonnegative constant\\
\begin{equation}
\alpha_2=\sum_{q\geq 1} \left | \frac{\mu(q)}{q^{s-1}\varphi(q)} \right |^2c_{q}(2)\geq \frac{1}{2}
\end{equation} \\
for $\sigma>1$ shows that the subset of quadratic twin primes $\{ n^{2}+1, n^{2}+3: n\geq1 \}$ has positive density in the set of quadratic primes $\mathbb{P}_{2} = \{ 5, 17, 37,\ldots, n^{2}+1, \ldots :n\geq1\}$.
\\

\subsection{Numerical Data}
A small numerical experiment was conducted to show the abundance of quadratic twin primes and to estimate the constant. The counting function is defined by\\
\begin{equation}
\pi _{f_2}(x)=\#\{p=n^2+1\leq x: p \text{ and } p+2 \text{ are primes}\}.
\end{equation}\\
The data shows a converging stable constant. \\

\begin{center}
 \begin{tabular}{||c |c| c| c|c|c|c|} 
 \hline
$ x$ & $\pi_{f_2}(x)$  & $\pi_{f_2}(x)/(x^{1/2}/\log x)$ & $ x$ & $\pi_{f_2}(x)$  & $\pi_{f_2}(x)/(x^{1/2}/\log x)$ \\ [1ex] 
 \hline\hline
$ 10^3$ & 4 &6.03579 &$ 10^9$&278&3.77538\\ 
 \hline
$10^5$&12 &5.02982 &$ 10^{10}$&689 &3.65301 \\
 \hline
 $10^6$ &28 &5.34431 &$ 10^{11}$&1782 &3.61513 \\
 \hline
$10^7$ & 6 1& 5.01138 &$ 10^{12}$&4663 &3.56008\\
 \hline
 $10^8$ &120 &4.07186 & $ 10^{13}$& 12260&3.47383 \\ 
 \hline
\end{tabular}
\end{center}

\vskip 1 in
\section{Cubic Twin Primes}
The cubic twin primes conjecture calls for an infinite subset of prime pairs $n^3+2, n^3+4$ as the integer $n \geq 1$ tends to infinity. There is no information on this primes counting problem available in the literature. \\

\textbf{Cubic Twin Primes Conjecture.}  There are infinitely many twin primes $n^3+2, n^3+4$ as $n \to \infty$. Moreover, the counting function has the asymptotic formula\\
\begin{equation}
\pi _{f_3}(x)=2\frac{x^{1/3}}{\log ^2 x}\prod _{p\geq 5} \left(\frac{p(p-\varrho (p))}{(p-1)^2}\right)+O\left(\frac{x^{1/3}}{\log ^3 x}\right),
\end{equation}\\
where \(\varrho (p)=0,2,6\) according to the occurrence of $0,1,\text{ or }2$ cubic residues $\text{mod } p$ in the set $\{-2,-4\}$.\\

The heuristic based on probability, yields the Gaussian type analytic formula\\
\begin{equation}
\pi _{f_3}(x)=s_{3} \int_{2}^{x} \frac{1}{ t^{2/3} \log^2 t} dt+O\left(\frac{x^{1/3}}{\log ^3 x}\right),
\end{equation}\\
where $s_{3}$ is the same constant as above. The deterministic approach is based on the weighted prime pairs counting function\\
\begin{equation}
\sum_{n\leq x}\Lambda(n^3+2)\Lambda(n^3+4),
\end{equation}\\
which is basically an extended version of the Chebishev method for counting primes.\\

\subsection{The Proof of the Theorem}
A lower bound for the correlation function $\sum_{n^{3}+2\leq x}\Lambda(n^{3}+2)\Lambda(n^{3}+4)$ will be derived here. It exploits the duality or equivalence of the Ramanujan series domain and the Taylor series domain to compute the correlation function of in two distinct ways. \\

\textbf{Theorem 10.2.}  \textit{Let $x\geq 1$ be a large real number. Then}\\
\begin{equation}
\sum_{n^{3}+2\leq x}\Lambda(n^{3}+2)\Lambda(n^{3}+4) \geq  \alpha_{2} x^{1/3}+O\left (x^{1/3} \frac{(\log \log x)^2}{\log x} \right ),
\end{equation}\\
\textit{for some $\alpha_{2} > 1/2$ constant.}  \\

\textbf{Proof:} Let $x \geq 1$ be a large number, and let $\varepsilon= ({\log \log x)^2}/{\log x}>0$. By Lemma 2.1, the function $f(n^{3}+2,s)$ has an absolutely, and uniformly convergent Ramanujan series for 
$\Re e(s) = \sigma \geq 1+c({\log \log x)^2}/{\log x}>1$, with $0<c<1$ constant. \\

By Lemma 6.1, the correlation function of the Taylor series of $f(n^{3}+2,s)$ has the asymptotic formula\\
\begin{equation} \label{el390}
\sum_{n^{3}+2\leq x} f(n^{3}+1,s)f(n^{3}+4,s) 
= \sum_{n^{3}+2\leq x}\Lambda(n^{3}+2)\Lambda(n^{3}+4)+O \left (x^{1/3}\frac{(\log\log x)^2}{\log x} \right ),
\end{equation}\\
where $0<|s-1|<\delta$, with $0< \delta \ll \varepsilon=({\log \log x)^2}/{\log x}$, see Lemma 4.4.\\

Applying Lemma 6.2, to the correlation function of the Ramanujan series of $f(n^{3}+2,s)$ returns the asymptotic formula\\
\begin{equation} \label{el391}
\sum_{n^{3}+2\leq x} f(n^{3}+2,s)f(n^{3}+4,s) 
=\alpha_{2} x^{1/3}+O \left (x^{1/3}\frac{(\log\log x)^2}{\log x} \right ),
\end{equation}
\\
where $\Re e(s)=1+c(\log\log x)^{2}/\log x=1+2\beta>1$ as required by Theorem 6.2. \\ 

Since the constant $\alpha_k>1/2$ for $\Re e(s)>1$, see Theorem 7.1, combining (\ref{el390}) and (\ref{el391}) confirms the claim: \\
\begin{eqnarray}
\sum_{n^3+2\leq x}\Lambda(n^{3}+2)\Lambda(n^{3}+4) 
&=& \sum_{n^{3}+2\leq x} f(n^{3}+2,s)f(n^{3}+4,s)  \\
&=& \alpha_{2}x^{1/3}+O \left(x^{1/3}\frac{(\log\log x)^{2}}{\log x}\right)\\
&\geq& (1/2)x^{1/3}+O \left(x^{1/3}\frac{(\log\log x)^{2}}{\log x}\right),
\end{eqnarray}\\
for a lower bound of the correlation function.   $\blacksquare$  \\

The nonnegative constant\\
\begin{equation}
\alpha_2=\sum_{q\geq 1} \left | \frac{\mu(q)}{q^{s-1}\varphi(q)} \right |^2c_{q}(2)\geq \frac{1}{2}
\end{equation} \\
for $\sigma>1$ shows that the subset of cubic twin primes $\{ n^{3}+2, n^{3}+4: n\geq1 \}$ has positive density in the set of quadratic primes $\mathbb{P}_{3} = \{ 1^1+2, 3^3+2, 45^3+2,63^3+2\ldots, n^{3}+2, \ldots :n\geq1\}$.
\\

\subsection{Numerical Data}
A small numerical experiment was conducted to show the abundance of cubic twin primes and to estimate the constant. The counting function is defined by\\
\begin{equation}
\pi _{f_3}(x)=\#\{p=n^3+2\leq x: p \text{ and } p+2 \text{ are primes}\}.
\end{equation}\\
The data shows an unstable constant. Perhaps it is not sufficient to determine a good approximation.\\

\begin{center}
 \begin{tabular}{||c |c| c| c|c|c|c|} 
 \hline
$ x$ & $\pi_{f_3}(x)$  & $\pi_{f_3}(x)/(x^{1/3}/\log x)$ & $ x$ & $\pi_{f_3}(x)$  & $\pi_{f_3}(x)/(x^{1/3}/\log x)$ \\ [1ex] 
 \hline\hline
$ 10^3$ & 2 & 9.54342 &$ 10^9$&12&5.15344 \\ 
 \hline
$10^5$&3 &8.56694 &$ 10^{10}$&17 &4.18357 \\
 \hline
 $10^6$ &5 &9.54342 &$ 10^{11}$&30 &4.14640 \\
 \hline
$10^7$ & 6 & 7.23511 &$ 10^{12}$&49 &3.74102\\
 \hline
 $10^8$ &7 &5.11732 & $ 10^{13}$& 83&3.45194 \\ 
 \hline
\end{tabular}
\end{center}

\vskip 1 in
\section{Quartic Twin Primes}
The quartic twin primes conjecture calls for an infinite subset of prime pairs $n^4+1, n^4+3$ as the integer $n \geq 1$ tends to infinity. There is no information on this primes counting problem available in the literature. \\

\textbf{Quartic Twin Primes Conjecture.}  There are infinitely many twin primes $n^4+1, n^4+3$ as $n \to \infty$. Moreover, the counting function has the asymptotic formula\\
\begin{equation}
\pi _{f_4}(x)=\frac{x^{1/4}}{\log ^2 x}\prod _{p\geq 2} \left(\frac{p(p-\varrho (p))}{(p-1)^2}\right)+O\left(\frac{x^{1/4}}{\log ^3 x}\right),
\end{equation}\\
where \(\varrho (p)=0,2,4\) according to the occurrence of $0,1,\text{ or }2$ quartic residues $\text{mod } p$ in the set $\{-1,-3\}$.\\

The heuristic based on probability, yields the Gaussian type analytic formula\\
\begin{equation}
\pi _{f_4}(x)=t_{4} \int_{2}^{x} \frac{1}{ t^{3/4} \log^2 t} dt+O\left(\frac{x^{1/4}}{\log ^3 x}\right),
\end{equation}\\
where $t_{4}$ is the same constant as above. The deterministic approach is based on the weighted prime pairs counting function\\
\begin{equation}
\sum_{n\leq x}\Lambda(n^4+1)\Lambda(n^4+3),
\end{equation}\\
which is basically an extended version of the Chebishev method for counting primes.\\

\subsection{The Proof of the Theorem}
A lower bound for the correlation function $\sum_{n^{4}+1\leq x}\Lambda(n^{4}+1)\Lambda(n^{4}+3)$ will be derived here. It exploits the duality or equivalence of the Ramanujan series domain and the Taylor series domain to compute the correlation function of in two distinct ways. \\

\textbf{Theorem 11.2.}  \textit{Let $x\geq 1$ be a large real number. Then}\\
\begin{equation}
\sum_{n^{4}+1\leq x}\Lambda(n^{4}+1)\Lambda(n^{4}+3) \geq  \alpha_{2} x^{1/4}+O\left (x^{1/4} \frac{(\log \log x)^2}{\log x} \right ),
\end{equation}\\
\textit{for some $\alpha_{2} > 1/2$ constant.}      \\

\textbf{Proof:} Let $x \geq 1$ be a large number, and let $\varepsilon= ({\log \log x)^2}/{\log x}>0$. By Lemma 2.1, the function $f(n^{4}+1,s)$ has an absolutely, and uniformly convergent Ramanujan series for 
$\Re e(s) = \sigma \geq 1+c({\log \log x)^2}/{\log x}>1$, with $0<c<1$ constant. \\

By Lemma 6.1, the correlation function of the Taylor series of $f(n^{3}+2,s)$ has the asymptotic formula\\
\begin{equation} \label{el490}
\sum_{n^{4}+1\leq x} f(n^{4}+1,s)f(n^{4}+3,s) 
= \sum_{n^{4}+1\leq x}\Lambda(n^{4}+1)\Lambda(n^{4}+3)+O \left (x^{1/4}\frac{(\log\log x)^2}{\log x} \right ),
\end{equation}\\
where $0<|s-1|<\delta$, with $0< \delta \ll \varepsilon=({\log \log x)^2}/{\log x}$, see Lemma 4.4.\\

Applying Lemma 6.2, to the correlation function of the Ramanujan series of $f(n^{4}+1,s)$ returns the asymptotic formula\\
\begin{equation} \label{el491}
\sum_{n^{4}+1\leq x} f(n^{4}+1,s)f(n^{4}+3,s) 
=\alpha_{2} x^{1/4}+O \left (x^{1/4}\frac{(\log\log x)^2}{\log x} \right ),
\end{equation}
\\
where $\Re e(s)=1+c(\log\log x)^{2}/\log x=1+2\beta>1$ as required by Theorem 6.2. \\ 

Since the constant $\alpha_k>1/2$ for $\Re e(s)>1$, see Theorem 7.1, combining (\ref{el490}) and (\ref{el491}) confirms the claim: \\
\begin{eqnarray}
\sum_{n^4+1\leq x}\Lambda(n^{4}+1)\Lambda(n^{4}+3) 
&=& \sum_{n^{4}+1\leq x} f(n^{4}+1,s)f(n^{4}+3,s)  \\
&=& \alpha_{2}x^{1/4}+O \left(x^{1/4}\frac{(\log\log x)^{2}}{\log x}\right)\\
&\geq& (1/2)x^{1/4}+O \left(x^{1/4}\frac{(\log\log x)^{2}}{\log x}\right),
\end{eqnarray}\\
for a lower bound of the correlation function.   $\blacksquare$  \\

The nonnegative constant\\
\begin{equation}
\alpha_2=\sum_{q\geq 1} \left | \frac{\mu(q)}{q^{s-1}\varphi(q)} \right |^2c_{q}(2)\geq \frac{1}{2}
\end{equation} \\
for $\sigma>1$ shows that the subset of cubic twin primes $\{ n^{4}+1, n^{4}+3: n\geq1 \}$ has positive density in the set of quadratic primes $\mathbb{P}_{f_4} = \{ 2^4+1, 2^4+3, 16^4+1,16^4+3, 28^4+1,28^4+3,\ldots, n^{3}+2, \ldots :n\geq1\}$.
\\

\subsection{Numerical Data}
A small numerical experiment was conducted to show the abundance of cubic twin primes and to estimate the constant. The counting function is defined by\\
\begin{equation}
\pi _{4}(x)=\#\{p=n^4+1\leq x: p \text{ and } p+2 \text{ are primes}\}.
\end{equation}\\
The data shows an unstable constant. Perhaps it is not sufficient to determine a good approximation.\\

\begin{center}
	\begin{tabular}{||c |c| c| c|c|c|c|} 
		\hline
		$ x$ & $\pi_{f_4}(x)$  & $\pi_{f_4}(x)/(x^{1/4}/\log x)$ & $ x$ & $\pi_{f_4}(x)$  & $\pi_{f_4}(x)/(x^{1/4}/\log x)$ \\ [1ex] 
		\hline\hline
		$ 10^3$ & 1 & 8.48543 &$ 10^12$&11&8.39821 \\ 
		\hline
		$10^5$&2 &14.9074 &$ 10^{13}$&17 &8.56578 \\
		\hline
		$10^6$ &3 &18.1074 &$10^{14}$&28 &9.20122 \\
		\hline
		$10^7$ &4 &18.4794 &$ 10^{15}$&40 &8.48543\\
		\hline
		$10^8$ &5 &16.9661 & $ 10^{16}$& 63&8.5509 \\ 
		\hline
	\end{tabular}
\end{center}

\vskip 1 in
\section{Primes Values Of Nonlinear Polynomials}

\textbf{Definition 12.1.} The \textit{fixed divisor} $\text{div}(f)=\gcd(f(\mathbb{Z}))$ of a polynomial $f(x) \in \mathbb{Z}[x]$ over the integers is the greatest common divisor of its image $f(\mathbb{Z})=\{f(n):n \in \mathbb{Z}\}$. \\

The fixed divisor $\text{div}(f)=1$ if $f(x) \equiv 0 \mod p$ has $\rho(p)<p$ solutions for every prime $p<\deg(f)$, see \cite[p. 395]{FI10}. An irreducible polynomial can represent infinitely many primes if and only if it has fixed divisor $\text{div}(f)=1$.\\

\textbf{Example 12.2.} The polynomial $f_1(x)=x^2+1$ is irreducible over the integers and has the fixed divisor $\text{div}(f_1)=1$, so it can represent infinitely many primes. On the other hand, $f_2(x)=x(x+1)+2$ is irreducible over the integers and has the fixed divisor $\text{div}(f_2)=2$, so it cannot represents infinitely many primes.\\

\vskip 1 in
\section{Quadratic Polynomials}
The heuristic analysis for the expected number of primes $p=n^2+1 \leq x$ is based on the circle method. This derivation appears in \cite[p. 48]{HL23}. Introductions to the circle methods, at various technical levels, are given in \cite[Chapter 23]{MT06}, and \cite{VR97}.\\

\textbf{Conjecture 13.1.} (\cite{HL23}) The expected number of quadratic primes $p=n^2+1 \leq x$ has the asymptotic weighted counting function\\
\begin{equation}
\sum_{n^2+1\leq x}\Lambda(n^{2}+1) \geq C_f x^{1/2}+O \left(x^{1/2}/\log x\right),
\end{equation}\\
where the constant, in terms of the quadratic symbol $(-1|p)=(-1)^{(p-1)/2}$ is\\  
\begin{equation}
C_f=\prod_{p>2} \left (1-\dfrac{(-1|p)}{p-1} \right )=1.37281346 \ldots.
\end{equation}
This conjecture, also known as one of the Landau primes counting problems, and other related problems are discussed in \cite{BH62}, \cite{BR07}, \cite{GM00}, \cite{GL10}, \cite{MA09}, \cite[p. 343]{NW00}, \cite[p. 405]{RP96}, \cite[p. 13]{CC00}, \cite{PJ09}, etc.\\

\textbf{Theorem 13.2.}  There are infinitely many quadratic primes $p=n^2+1$ as $n \to \infty$. Moreover, its counting function has the weighted asymptotic formula\\
\begin{equation}
\sum_{n^2+1\leq x}\Lambda(n^{2}+1) \geq \alpha_{2}\frac{x^{1/2}}{\log x}+O \left(x^{1/2}\frac{(\log\log x)^{2}}{\log^2 x}\right),
\end{equation}\\
where $\alpha_2>1/2$ is a constant, for all large numbers $x \geq 1$.\\

\textbf{Proof:} Let $R(x)=\sum_{n^2+1\leq x}\Lambda(n^{2}+1)\Lambda(n^{2}+3)$. Since the weighted counting measure $R(x)$ under counts the primes $p=n^2+1\leq x$ by a factor of $\log x$, it can be used to estimate a lower bound for the weighted counting measure $\sum_{n^2+1\leq x}\Lambda(n^{2}+1)$ of the primes $p=n^2+1\leq x$. Specifically, by partial summation, this is\\
\begin{eqnarray}
\sum_{n^2+1\leq x}\Lambda(n^{2}+1)
&\geq&\sum_{n^2+1\leq x}\Lambda(n^{2}+1)\frac{\Lambda(n^{2}+3)}{\log(n^{2}+3)}\\
&=& \int_1^x \frac{1}{\log(t^{2}+3)}d R(t).
\end{eqnarray}\\
Using Theorem 1.2 to evaluate the Stieltjes integral returns\\
\begin{equation}
\sum_{n^2+1\leq x}\Lambda(n^{2}+1)
\geq \alpha_2 \frac{x^{1/2}}{\log x}+
O\left (x^{1/2} \frac{(\log \log x)^2}{\log x} \right ),
\end{equation}\\
where $\alpha_2>1/2$ is a constant, and $x\geq1$ is a large number.   $\blacksquare$.\\

Some numerical data and experiments are reported in \cite[p. 50]{HL23}, \cite{SD60}, \cite{CC00}, and \cite{RI15}.

\vskip 1 in
\section{Quartic Polynomials}
The heuristic analysis for the expected number of primes $p=n^4+1 \leq x$ is based on the Bateman-Horn conjecture, described in the penultimate Section. \\

This heuristic, which is based on probability, yields the Gaussian type analytic formula\\
\begin{equation}
\pi _{f_4}(x)=s_{4} \int_{2}^{x} \frac{1}{ t^{3/4} \log^2 t} dt+O\left(\frac{x^{1/4}}{\log ^3 x}\right),
\end{equation}\\
where $s_{4}$ is a constant. The deterministic approach is based on the weighted prime pairs counting function\\
\begin{equation}
\sum_{n\leq x}\Lambda(n^4+1)\Lambda(n^4+3),
\end{equation}\\
which is basically an extended version of the Chebishev method for counting primes.\\

\textbf{Conjecture 14.1.} The expected number of quartic primes $p=n^4+1 \leq x$ has the asymptotic weighted counting function\\
\begin{equation}
\sum_{n^4+1\leq x}\Lambda(n^{4}+1) = s_4 x^{1/4}+O \left(\frac{x^{1/4}}{\log^2 x}\right),
\end{equation}\\
where the constant, in terms of a few values of the quadratic symbols $(a|p)\equiv a^{(p-1)/2} \text{ mod } p$, is\\  
\begin{equation}
s_4=\frac{1}{4}\prod_{p>2} \left (1-\dfrac{(-1|p)+(-2|p)+(2|p)}{p-1} \right )=.66974 \ldots.
\end{equation}\\
Numerical approximation for this constant are given in \cite{SD60}.\\

This conjecture is discussed in \cite{SD60}, etc.\\

\subsection{Lower Bounds}
Lower bounds for the counting functions are computed here. These estimates have the correct order of magnitude expected up to a $\log x$ factor.\\
 
\textbf{Theorem 14.2.}  There are infinitely many quartic primes $p=n^4+1$ as $n \to \infty$. Moreover, its counting function has the weighted asymptotic formula\\
\begin{equation}
\sum_{n^4+1\leq x}\Lambda(n^{4}+1) \geq \alpha_{2}\frac{x^{1/4}}{\log x}+O \left(x^{1/4}\frac{(\log\log x)^{2}}{\log^2 x}\right),
\end{equation}\\
where $\alpha_2>1/2$ is a constant, for all large numbers $x \geq 1$.\\

\textbf{Proof:} Let $R(x)=\sum_{n^4+1\leq x}\Lambda(n^{4}+1)\Lambda(n^{4}+3)$. Since the weighted counting measure $R(x)$ under counts the primes $p=n^4+1\leq x$ by a factor of $\log x$, it can be used to estimate a lower bound for the weighted counting measure $\sum_{n^4+1\leq x}\Lambda(n^{4}+1)$ of the primes $p=n^4+1\leq x$. Specifically, by partial summation, this is\\
\begin{eqnarray}
\sum_{n^4+1\leq x}\Lambda(n^{4}+1)
&\geq&\sum_{n^4+1\leq x}\Lambda(n^{4}+1)\frac{\Lambda(n^{4}+3)}{\log(n^{4}+3)}\\
&=& \int_1^x \frac{1}{\log(t^{4}+3)}d R(t).
\end{eqnarray}\\
Using Theorem 1.4 to evaluate the Stieltjes integral returns\\
\begin{equation}
\sum_{n^4+1\leq x}\Lambda(n^{4}+1)
\geq \alpha_2 \frac{x^{1/4}}{\log x}+
O\left (x^{1/4} \frac{(\log \log x)^2}{\log x} \right ),
\end{equation}\\
where $\alpha_2>1/2$ is a constant, and $x\geq1$ is a large number.   $\blacksquare$.\\

Some numerical data and experiments are reported in \cite{SD60}, and \cite{RI15}.


\newpage
 \begin {thebibliography}{999}

\bibitem{AP76} Apostol, Tom M. Introduction to analytic number theory. Undergraduate Texts in Mathematics. Springer-Verlag, New York-Heidelberg, 1976.

\bibitem{BH62} Bateman, P. T., Horn, R. A. (1962), A heuristic asymptotic formula concerning the distribution of prime numbers, Math. Comp. 16 363-367.

\bibitem{BR07} Baier, Stephan; Zhao, Liangyi. Primes in quadratic progressions on average. Math. Ann. 338 (2007), no. 4, 963-982. 

\bibitem{CC00} Chris K. Caldwell, An Amazing Prime Heuristic, Preprint, 2000.

\bibitem{CM15} Giovanni Coppola, M. Ram Murty, Biswajyoti Saha, On the error term in a Parseval type formula in the theory of Ramanujan expansions II, arXiv:1507.07862.

\bibitem{CO07} Cohen, Henri. Number theory. Vol. II. Analytic and modern tools. Graduate Texts in Mathematics, 240. Springer, New York, 2007.

\bibitem{DE76} Delange, Hubert. On Ramanujan expansions of certain arithmetical functions. Acta Arith. 31 (1976), no. 3, 259-270.

\bibitem{DH00} Davenport, Harold Multiplicative number theory. Third edition. Revised and with a preface by Hugh L. Montgomery. Graduate Texts in Mathematics, 74. Springer-Verlag, New York, 2000.

\bibitem{DLMF} Digital Library Mathematical Functions, \text{http://dlmf.nist.gov.}

\bibitem{EL85}  Ellison, William; Ellison, Fern. Prime numbers. A Wiley-Interscience Publication. John Wiley and Sons, Inc., New York; Hermann, Paris, 1985.

\bibitem{FI10} Friedlander, John; Iwaniec, Henryk. Opera de cribro. American Mathematical Society Colloquium Publications, 57. American Mathematical Society, Providence, RI, 2010.

 \bibitem{FK12} Christopher F. Fowler, Stephan Ramon Garcia, Gizem Karaali, Ramanujan sums as supercharacters, arXiv:1201.1060.

\bibitem{GD09} D. A. Goldston, Are There Infinitely Many Twin Primes?, Preprint, 2009.

\bibitem{GG05} Goldston, D. A.; Graham, S. W.; Pintz, J.; Yildirim, C. Y. Small gaps between almost primes, the parity problem, and some conjectures of Erdos on consecutive integers, arXiv:math.NT/0506067.

\bibitem{GL10} Luca Goldoni, Prime Numbers And Polynomials, Phd Thesis, Universita` Degli Studi Di Trento, 2010.

\bibitem{GM00} Granville, Andrew; Mollin, Richard A. Rabinowitsch revisited. Acta Arith. 96 (2000), no. 2, 139-153. 

\bibitem{GP99} Gadiyar, H. Gopalkrishna; Padma, R. Ramanujan-Fourier series, the Wiener-Khintchine formula and the distribution of prime pairs. Phys. A  269  (1999),  no. 2-4, 503-510.

\bibitem{GP06} H. Gopalkrishna Gadiyar, R. Padma, Linking the Circle and the Sieve: Ramanujan - Fourier Series, arXiv:math/0601574.

\bibitem{HG07} Harman, Glyn. Prime-detecting sieves. London Mathematical Society Monographs Series, 33. Princeton University Press, Princeton, NJ, 2007.

\bibitem{HL23} Hardy, G. H. Littlewood J.E. Some problems of Partitio numerorum III: On the expression of a number as a sum of primes. Acta Math. 44 (1923), No. 1, 1-70.

\bibitem{HR59} Hardy, G. H. Ramanujan: twelve lectures on subjects suggested by his life and work.  Chelsea Publishing Company, New York 1959.

\bibitem{HW08} Hardy, G. H.; Wright, E. M. An introduction to the theory of numbers. Sixth edition. Revised by D. R. Heath-Brown and J. H. Silverman. With a foreword by Andrew Wiles. Oxford University Press, Oxford, 2008.

\bibitem{JW03} Jacobson, Michael J., Jr.; Williams, Hugh C. New quadratic polynomials with high densities of prime values. Math. Comp. 72 (2003), no. 241, 499-519.

\bibitem{KJ92} Keiper, J. B. Power series expansions of Riemann's $\xi$ function. Math. Comp. 58 (1992), no. 198, 765-773.

\bibitem{KL12} Kunik, Matthias; Lucht, Lutz G. Power series with the von Mangoldt function. Funct. Approx. Comment. Math. 47 (2012),  part 1, 15-33. 

\bibitem{KR09} Jacob Korevaar. Prime pairs and the zeta function, Journal of Approximation Theory 158 (2009) 69-96.

\bibitem{LL95} Lucht, Lutz Ramanujan expansions revisited. Arch. Math. (Basel) 64 (1995), no. 2, 121-128.

\bibitem{LL10} Lucht, Lutz G. A survey of Ramanujan expansions. Int. J. Number Theory 6 (2010), no. 8, 1785-1799. 

\bibitem{MA09} Matomaki, Kaisa. A note on primes of the form p=aq2+1. Acta Arith. 137 (2009), no. 2, 133-137.

\bibitem{MT06} Miller, Steven J.; Takloo-Bighash, Ramin. An invitation to modern number theory. With a foreword by Peter Sarnak. Princeton University Press, Princeton, NJ, 2006. 

\bibitem{MV07} Montgomery, Hugh L.; Vaughan, Robert C. Multiplicative number theory. I. Classical theory. Cambridge University Press, Cambridge, 2007.

\bibitem{NW00} Narkiewicz, W. The development of prime number theory. From Euclid to Hardy and Littlewood. Springer Monographs in Mathematics. Springer-Verlag, Berlin, 2000. 

\bibitem{PJ09} Pintz, Janos. Landau's problems on primes. J. Theory. Nombres Bordeaux 21 (2009), no. 2, 357-404.

\bibitem{RM13} Ram Murty, M. Ramanujan series for arithmetical functions. Hardy-Ramanujan J.  36  (2013), 21-33.

\bibitem{RP96} Ribenboim, Paulo, The new book of prime number records, Berlin, New York: Springer-Verlag, 1996.

\bibitem{RI15} Igor Rivin, Some experiments on Bateman-Horn,  arXiv:1508.07821.

\bibitem{SD60} Shanks, Daniel On numbers of the form $n^4 +1$ . Math. Comput.  15,  1961, 186-189.

\bibitem{SD61} Shanks, Daniel A sieve method for factoring numbers of the form $n^2 +1$ . Math. Tables Aids Comput.  13,  1959, 78-86.

\bibitem{SP80} Spilker, Jurgen. Ramanujan expansions of bounded arithmetic functions. Arch. Math. (Basel) 35 (1980), no. 5, 451-453.

\bibitem{SS94} Schwarz, Wolfgang; Spilker, Jurgen. Arithmetical functions. An introduction to elementary and analytic properties of arithmetic functionss. London Mathematical Society Lecture Note Series, 184. Cambridge University Press, Cambridge, 1994.

\bibitem{TG15} Tenenbaum, Gerald. Introduction to analytic and probabilistic number theory. Translated from the Third French edition. American Mathematical Society, Rhode Island, 2015.

\bibitem{VR97} Vaughan, R. C. The Hardy-Littlewood method. Second edition. Cambridge Tracts in Mathematics, 125. Math. Soc. (N.S.) 44 (2007), no. 1, 1-18. Cambridge University Press, Cambridge, 1997. 

\bibitem{WN66} Wiener, Norbert. Generalized harmonic analysis and Tauberian theorems. M.I.T. Press, Cambridge, Mass., London 1966.



\end{thebibliography}
\end{document}